\documentclass[11pt,twoside]{article}
\usepackage{amsmath}
\usepackage{amsthm}
\usepackage{amssymb}
\usepackage{array}
\usepackage{geometry}
\usepackage{graphicx}
\usepackage{enumerate}
\usepackage{lscape,graphicx}
\usepackage{fancyhdr}
\usepackage{makeidx}
\usepackage{amsfonts}
\usepackage{graphicx}
\usepackage{CJK}
\usepackage{amsmath}
\usepackage{latexsym}
\usepackage{amssymb}
\usepackage{mathrsfs}
\usepackage{amsthm}
\usepackage{array}
\usepackage{enumerate}
\usepackage{longtable}
\usepackage{multirow}
\makeindex
\usepackage{wasysym}
\usepackage{chapterbib}

\normalsize
\begin{document}


\newpage
\setcounter{equation}{0}
\begin{center}
\vskip1cm
{\Large
\textbf{Some Properties and Applications of Burr III-Weibull Distribution} }
\vskip.5cm
%

\vspace{0.50cm}
\textbf{\large Deepthy G S, Nicy Sebastian  and Reshma Rison }\\
Department of Statistics,
St.Thomas College, Thrissur, Kerala, India\\ {Email:{\tt{deepthygs@gmail.com, nicycms@gmail.com}}}\\

\end{center}

\thispagestyle{empty}
\begin{center} {\bf }  \vskip 0.10truecm \noindent \\

\vskip.5cm\noindent \centerline{\bf Abstract}
\end{center}
\par
In this paper, we introduce a new distribution called Burr III-Weibull(BW) distribution using the concept of competing risk. We derive  moments, conditional moments, mean deviation and quantiles of the proposed distribution. Also the Renyi's entropy and order statistics of  the distribution are obtained. Estimation of parameters
of the distribution is performed via maximum likelihood method. A simulation study is performed to validate the maximum likelihood estimator (MLE). A real practical data set
is analyzed for illustration. 		 \\



\vskip.5cm{\section*{1.\hskip.3cm Introduction}}
\vskip.3cm

Burr type III distribution with two parameters was first introduced in the literature of Burr \cite{1} for modelling lifetime data or survival data. It is more flexible and includes a variety of distributions with varying degrees of skewness and kurtosis. This distribution has a wide application in areas of statistical modelling such as forestry Gove et al.\cite{3}, meteorology Mielke \cite{2}, and reliability Mokhlis \cite{7}. Burr type XII distribution can be derived from Burr type III distribution by replacing $X$ with $\frac{1}{X}$. The usefulness and properties of Burr distribution are discussed by Burr and Cislak \cite{8} and Johnson et al. \cite{9}.

	\vskip.3cm The hazard function of distributions may include one or more of the following behavioural patterns; increasing, decreasing or constant
	shapes. Thus, they cannot be used to model lifetime data with a bathtub shaped hazard function, such as human mortality and machine life cycles. For last few decades, statisticians have been developing various extensions and modified forms of the Weibull distribution and other related models.
	The two-parameter, flexible Weibull extension of Bebbington et al. \cite{4}  has a hazard function that can be increasing, decreasing or bathtub shaped. Zhang and Xie \cite{5}  studied the characteristics
	and application of the truncated Weibull distribution, which has a bathtub shaped hazard function. A new modified Weibull distribution by Saad and Jingsong \cite{6}  considered an increasing and a bathtub shaped hazard function.
\vskip.3cm
		The cumulative distribution function(cdf) and probability density function(pdf) of the weibull distribution are given by,
	\begin{eqnarray}
	F_{W}(x;\lambda,\beta)&=&1-e^{-\left(x/\lambda\right)^{\beta}}, x\geq0,\lambda>0,\beta>0
	\\f_W(x;\lambda,\beta) &=& \frac{\beta}{\lambda}\left(\frac{x}{\lambda}\right)^{\beta-1}e^{-\left(x/\lambda\right)^{\beta}}
	\end{eqnarray}
where $\lambda$ and $\beta$ are the scale and shape parameters. The cumulative distribution(cdf) and probability density function(pdf) of the BurrIII distribution is given by,
	\begin{eqnarray}G_B(x;c,k)&=&(1+x^{-c})^{-k}, x\geq0, k>0, c>0
	\\g_B(x;c,k)&=& ckx^{-c-1}\left(1+x^{-c}\right)^{-k-1}
	 \end{eqnarray} where c and k are shape parameters.

\vskip.5cm{\section*{2.\hskip.3cm Burr III-Weibull Distribution}}
\vskip.3cm
	
The reliability function of the new distribution, say Burr III Weibull(BW) distribution, can be constructed by combining the corresponding reliability functions of Burr III and Weibull distributions. The resulting reliability function, the cumulative distribution function and the probability density function are given by,
\begin{eqnarray}
\bar{F}_{BW}(x;c,k,\lambda,\beta) &=&\left(1-\left(1+x^{-c}\right)^{-k}\right)\left(e^{-\left(x/\lambda\right)^{\beta}}\right)  \ \ ; c,k,\lambda,\beta>0\\
F_{BW}(x;c,k,\lambda,\beta)&=&1-\left(1-(1+x^{-c})^{-k}\right)\left(e^{-\left(x/\lambda\right)^{\beta}}\right)
\end{eqnarray}for \text{c,k,$\lambda$,$\beta>0$}.
\begin{eqnarray}
	f_{BW}(x;c,k,\lambda,\beta)&=&e^{-\left(x/\lambda\right)^{\beta}} \left[ck\left(1+x^{-c}\right)^{-k-1}x^{-c-1}+\frac{\beta}{\lambda^\beta}x^{\beta-1}\left(1-\left(1+x^{-c}\right)^{-k}\right)\right]
\end{eqnarray}for $c,k,\lambda,\beta>0$.
The hazard rate \text{h(x)} and reverse hazard rate \text{r(x)} are given respectively as,
\begin{eqnarray}
h(x)&=&\frac{ck(1+x^{-c})^{-k-1}x^{-c-1}+\frac{\beta}{\lambda^{\beta}}x^{\beta-1}(1-(1+x^{-c})^{-k})}{1-(1+x^{-c})^{-k}}
\\T(x) &=&\frac{e^{-\left(x/\lambda\right)^{\beta}} [ck(1+x^{-c})^{-k-1}x^{-c-1}+\frac{\beta}{\lambda^{\beta}}x^{\beta-1}(1-(1+x^{-c})^{-k})]}{1-[1-(1+x^{-c})^{-k}e^{-\left(x/\lambda\right)^{\beta}}}
\end{eqnarray}

\begin{figure}
\begin{center}
	\includegraphics[height=.4\textheight,width=0.8\linewidth]{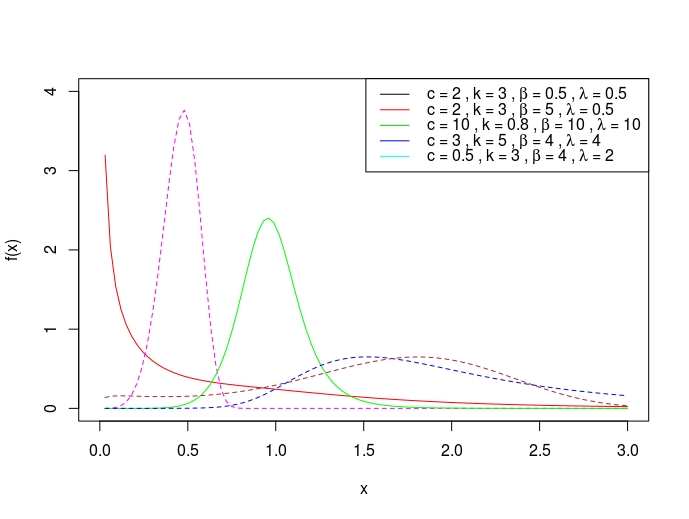}
	
	\caption{Plot for probability density functions of the BurrIII-Weibull distribution.}	
\end{center}
\end{figure}

\begin{figure}
\begin{center}
	\includegraphics[height=.4\textheight,width=0.8\linewidth]{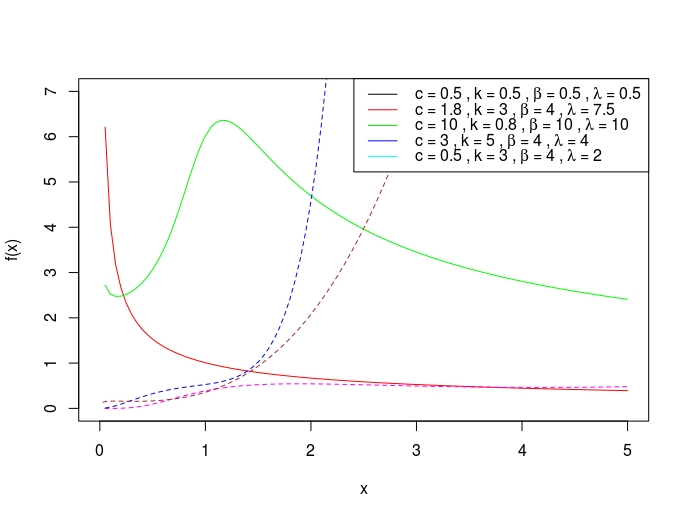}
	\caption{Plot for hazard rate functions of the BurrIII-Weibull distribution.}	
\end{center}
\end{figure}
The plots of the probability  and hazard rate functions of BW distribution for selected values of parameters are displayed in Figure 1 and Figure 2.  Figure 1 shows that the BW distribution can be decreasing, approximately symmetrical and right skewed whereas Figure 2 reflects the monotonic behaviour of the hazard function for different values of parameters.
%
%

\vskip.5cm{\section*{3.\hskip.3cm The Statistical Properties}}
\vskip.3cm
	
In this section, some of the statistical properties of BW distribution such as quantile function, moments and order statistics are derived.

\vskip.3cm {\subsection*{3.1\hskip.3cm Quantile Function}}\vskip.3cm

The quantile function has a number of important applications, for example, it can be used
to obtain the median, skewnes and kurtosis, and can also be used to generate random variables.
The quantile function can be obtained by inverting $\bar{F}_{BW}(x)=1-u, 0\leq u\leq 1$, where
\begin{eqnarray}\bar{F}_{BW}(x)&=&\left(1-(1+x^{-c})^{-k}\right)\left(e^{-\left(x/\lambda\right)^{\beta}}\right).  \end{eqnarray}
Let, $\left(1-(1+x^{-c})^{-k}\right)\left(e^{-\left(x/\lambda\right)^{\beta}}\right)=1-u$,
and the corresponding quantile function is obtained by solving the non-linear equation,
\begin{eqnarray}\label{quan}
\ln[1-(1+x^{-c})^{-k}]-\left(\frac{x}{\lambda}\right)^\beta-\ln(1-u)=0, \end{eqnarray}
using numerical methods. Equation (\ref{quan}) can be used to generate random number. The quantiles for selected
values of the BW distribution parameters are listed in Table \ref{t1}.
{\scriptsize{
\begin{center}
\begin{table}
	\caption{BW quantile for selected values}
	\label{t1}
	\begin{tabular}{c c c c c c}
	\hline
	&\multicolumn{5}{c}{(c,k,$\lambda$,$\beta$)}\\
	\hline	
	u&(3,1,2,0.4)&(0.1,1.7,1,1)&(1.8,1.3,0.6,3) &(3,0.1,0.7,0.5)&(0.5,1.2,1,0.8)\\
	\hline
	0.1&0.00720&0.00007&0.24762&0.00023& 0.014564\\
	0.2&0.04698&0.00776&0.32549&0.00202& 0.04755\\
	0.3&0.14851&0.07691&0.38578&0.00705& 0.09955\\
	0.4&0.31794&0.20024 &0.43923&0.01750&0.17508\\
	0.5&0.50091&0.36225 &0.49027&0.03678&0.28234\\
	0.6& 0.67537&0.56928&0.54192&0.070287&0.43596\\
	0.7& 0.85894&0.84254&0.59757&0.12815&0.66495\\
	0.8&1.08687&1.2336 &0.66279&0.23202& 1.03706\\
	0.9& 1.46385&1.90969&0.75279& 0.44542& 1.78102\\
	\end{tabular}
\end{table}
\end{center}}}

\vskip.3cm {\subsection*{3.2\hskip.3cm Moments}}\vskip.3cm

Moments can be used to study the most important features and characteristics of a distribution such as central tendency, dispersion, skewness, kurtusis etc.
The $r^{th}$ moment of BW distribution is given by,
\begin{eqnarray*}
E(X^r) &=& \int_{0}^{\infty}x^rf(x)dx
\\&=& \int_{0}^{\infty}x^{r}e^{-(x/\lambda)^{\beta}}\left[ck(1+x^{-c})^{-k-1}x^{-c-1}+\frac{\beta}{\lambda^\beta}x^{\beta-1}(1-(1+x^{-c})^{-k})\right]dx
\\&=&ck\int_{0}^{\infty}x^{r-c-1}(1+x^{-c})^{-k-1}e^{-({x/\lambda})^{\beta}}dx+\frac{\beta}{\lambda^{\beta}}\int_{0}^{\infty}x^{r+\beta-1}e^{-({x/\lambda})^{\beta}}dx\\&-& \frac{\beta}{\lambda^{\beta}}\int_{0}^{\infty}x^{r+\beta-1}(1+x^{-c})^{-k}e^{-({x/\lambda})^{\beta}}dx.
\end{eqnarray*}
Let
\begin{eqnarray*}
A&=&ck\int_{0}^{\infty}x^{r-c-1}(1+x^{-c})^{-k-1}e^{-({x/\lambda})^{\beta}}dx,\nonumber\\B&=&\frac{\beta}{\lambda^{\beta}}
\int_{0}^{\infty}x^{r+\beta-1}e^{-({x/\lambda})^{\beta}}dx \text{and} \nonumber\\
 C&=&\frac{\beta}{\lambda^{\beta}}\int_{0}^{\infty}x^{r+\beta-1}(1+x^{-c})^{-k}e^{-({x/\lambda})^{\beta}}dx.
\end{eqnarray*}
Then
\begin{equation} \label{m}
E(X^r)=A+B-C .
\end{equation}
Consider,
\begin{eqnarray}
	A&=&ck\sum_{m=0}^{\infty}\frac{(-1)^{m}}{\lambda^{m\beta}m!}\int_{0}^{\infty}\left(1+x^{-c}\right)^{-k-1}x^{r+m\beta-c-1}dx ,\text{put $u=(1+x^{-c})^{-1}$ ,}\nonumber\\ &=&ck\sum_{m=0}^{\infty}\frac{(-1)^{m}}{\lambda^{(m\beta)}m!}\int_{0}^{\infty}u^{k+\frac{r}{c}+\frac{m\beta}{c}-1}(1-u)^{1-\frac{r}{c}-\frac{m\beta}{c}-1}du \nonumber \\
	&=&ck\sum_{m=0}^{\infty}\frac{(-1)^{m}}{\lambda^{(m\beta)}m!}B\left(k+\frac{r}{c}+\frac{m\beta}{c},1-\frac{r}{c}-\frac{m\beta}{c}\right)\label{m1}
\end{eqnarray}
\begin{eqnarray}
B&=&\frac{\beta}{\lambda^{\beta}}\int_{0}^{\infty}x^{r+\beta-1}e^{-({x/\lambda})^{\beta}}dx\nonumber\\
&=&\frac{\beta}{\lambda^{\beta}}\frac{\lambda^{r+\beta}}{\beta}\Gamma \left(\frac{r+\beta}{\beta}\right),\text{using generalised gamma distribution}\label{m2} \\
C&=&\frac{\beta}{\lambda^{\beta}}\sum_{t=0}^{\infty}(-1)^{t}\binom{k+t-1}{t}\int_{0}^{\infty}x^{r+\beta-ct-1}e^{-(x/\lambda)^{\beta}}dx\nonumber \\
&=& \sum_{t=0}^{\infty}(-1)^{t}\binom{k+t-1}{t}\lambda^{r-ct}\Gamma\left(\frac{r+\beta-ct}{\beta}\right)\label{m3}
\end{eqnarray}
Substituting (\ref{m1}), (\ref{m2}) and (\ref{m3}) in (\ref{m}), we get,
\begin{eqnarray}
E(X^{r})&=&ck\sum_{m=0}^{\infty}\frac{(-1)^{m}}{\lambda^{(m\beta)}m!}B\left(k+\frac{r}{c}+\frac{m\beta}{c},1-\frac{r}{c}-\frac{m\beta}{c}\right)+\lambda^{r}\Gamma\left(\frac{r}{\beta}+1\right)\nonumber
\\&-&\sum_{t=0}^{\infty}(-1)^{t}\binom{k+t-1}{t}\Gamma\left(\frac{r+\beta-ct}{\beta}\right)\lambda^{r-ct}\ \text{where}\ \  r <c,m\beta<c.
\end{eqnarray}
Where $B\left(a,b\right)=\int_{0}^{1}t^{a-1}(1-t)^{b-1}dt$ is the beta function and $ \int_{0}^{\infty}x^{d-1}e^{-(x/a)^{p}}dx=\frac{\Gamma\left(\frac{d}{p}\right)}{\frac{p}{a^{d}}}$ is generalised gamma function. The moment generating function of the BW distribution is given by, $E\left(e^{tY}\right)= \sum_{i=0}^{\infty}\frac{t^{i}}{i!}E\left(Y^{i}\right)$ where $E\left(Y^{i}\right)$ is given above.
{\scriptsize{\begin{center}
\begin{table}
	\caption{BW moments for selected values $(c,k,\lambda,\beta)$}
	\label{moments}
	\begin{tabular}{c c c c c c}\hline
		Moments& (5,2.5,0.5,1.5)&(1,1,0.2,0.9)&(2,3,0.4,0.8)&(3,1.2,0.8,1.5)&(0.4,0.2,1,2)\\
		\hline
		$\mu_{1}^{'}$&0.44431&0.17564&0.21448&0.60806&0.14926\\
		$\mu_{2}^{'}$&0.28101&0.06958&0.12327&0.50286&0.14075\\
		$\mu_{3}^{'}$&0.21902&0.04434&0.12371&0.51011&0.17208\\
		$\mu_{4}^{'}$&0.19687&0.03965&0.18929&0.60844&0.24533\\
		$\mu_{5}^{'}$&0.19695&0.04618&0.41159&0.83411&0.39184\\
		$\mu_{6}^{'}$&0.21481&0.06674&1.20954&1.29465&0.68564\\
		SD& 0.28913&0.19680&0.27797&0.36486&0.34419\\
		CV&0.65073&1.12047&1.29601&0.60003&2.30597\\
		CS&0.82234&2.42891&2.98560&0.87391&2.83748\\
		Ck&3.37007&12.34781&18.5636&4.12879&11.39373\\
	\end{tabular}
\end{table}
\end{center}}}
The first six moments ($\mu_{1}^{'}$, $\mu_{2}^{'}$, $\mu_{3}^{'}$, $\mu_{4}^{'}$, $\mu_{5}^{'}$, $\mu_{6}^{'}$), standard deviation (SD), coefficient of variation (CV), coefficient of skewness (CS) and coefficient of kurtosis (CK) for different selected values of the BW distribution parameters are listed in Table \ref{moments}.
\vskip.3cm{\subsection*{3.3\hskip.3cm Conditional Moments}}

The $r^{th}$ conditional moment is defined as $E(X^r/X>t)$. The $r^{th}$ conditional moment of the BW distribution is given by,


\begin{eqnarray*}
&&E(X^{r}/X>t) \nonumber\\
&&\nonumber\\
& &\quad ~~~~=\frac{1}{\bar{F}(t)}\int_{t}^{\infty}x^{r}f(x)dx
\nonumber\\
& &\quad ~~~~=\frac{1}{\bar{F}(t)}\int_{t}^{\infty}x^{r}e^{-\left(\frac{x}{\lambda}\right)^{\beta}}
\left(ck\left(1+x^{-c}\right)^{-k-1}x^{-c-1}+\frac{\beta}{\lambda^{\beta}}x^{\left(\beta-1\right)}\left(1-\left(1+x^{-c}\right)^{-k}\right)\right)dx
\nonumber\\
& &\quad~~~~= \frac{1}{\bar{F}(t)}\int_{t}^{\infty}ck\left(1+x^{-c}\right)^{-k-1}x^{r-c-1}e^{-\left(\frac{x}{\lambda}\right)^{\beta}}dx+
\frac{1}{\bar{F}(t)}\frac{\beta}{\lambda^{\beta}}\int_{t}^{\infty}x^{r+\beta-1}e^{-\left(\frac{x}{\lambda}\right)^{\beta}}dx\nonumber\\
& &\quad~~~~ - \frac{1}{\bar{F}(t)}\frac{\beta}{\lambda^{\beta}}\int_{t}^{\infty}\left(1+x^{-c}\right)^{-k}e^{-\left(\frac{x}{\lambda}\right)^{\beta}}x^{r+\beta-1}dx.
\end{eqnarray*}
%

Let
\begin{eqnarray*}
  A&=&\int_{t}^{\infty}ck\left(1+x^{-c}\right)^{-k-1}x^{r-c-1}e^{-\left(\frac{x}{\lambda}\right)^{\beta}}dx,\nonumber\\
B&=&\frac{\beta}{\lambda^{\beta}}\int_{t}^{\infty}x^{r+\beta-1}e^{-\left(\frac{x}{\lambda}\right)^{\beta}}dx,\nonumber\\
C&=&\frac{\beta}{\lambda^{\beta}}\int_{t}^{\infty}\left(1+x^{-c}\right)^{-k}e^{-\left(\frac{x}{\lambda}\right)^{\beta}}x^{r+\beta-1}dx.\nonumber\\
\end{eqnarray*}
Then
\begin{equation}\label{c}
E(X^{r}/X>t)=\frac{1}{\bar{F}(t)}\left(A+B-C\right).
\end{equation}
Consider
\begin{eqnarray}\label{c1}
A&=& ck\sum_{m=0}^{\infty}\frac{(-1)^{m}}{\lambda^{m\beta}m!}\int_{t}^{\infty}\left(1+x^{-c}\right)^{-k-1}x^{r+m\beta-c-1}dx, \text{put}~~u=\left(1+x^{-c}\right)^{-1},\nonumber\\
&=&k\sum_{m=0}^{\infty}\frac{(-1)^{m}}{\lambda^{m\beta}m!}\int_{(1+t^{-c})^{-1}}^{1}u^{k+\frac{r}{c}+
  \frac{m\beta}{c}-1}\left(1-u\right)^{1-\frac{r}{c}-\frac{m\beta}{c}-1}du\nonumber\\
 &=&k\sum_{m=0}^{\infty}\frac{(-1)^{m}}{\lambda^{m\beta}m!}
 B_{(1+t^{-c})^{-1}}\left(k+\frac{r}{c}+\frac{m\beta}{c},1-\frac{r}{c}-\frac{m\beta}{c}\right)\label{c1}
\end{eqnarray}

%
\begin{eqnarray}\label{c2}
B&=&\frac{\beta}{\lambda^{\beta}}\int_{t}^{\infty}x^{r+\beta-1}e^{-\left(\frac{x}{\lambda}\right)^{\beta}}dx, \text{let}~~u=\left(\frac{x}{\lambda}\right)^{\beta}\nonumber\\
&=&\lambda^{r}\int_{(\frac{t}{\lambda})^{\beta}}^{\infty}u^{1+\frac{r}{\beta}-1}e^{-u}du\nonumber\\
&=&\lambda^{r}\Gamma\left(\left(\frac{r}{\beta}+1\right),\left(\frac{t}{\lambda}\right)^{\beta}\right)
\end{eqnarray}
\begin{eqnarray}\label{c3}
C&=&\frac{\beta}{\lambda^{\beta}}\sum_{p=0}^{\infty}(-1)^{p}\binom{k+p-1}{p}\int_{t}^{\infty}x^{r+\beta-cp-1}e^{-\left(\frac{x}{\lambda}\right)^{\beta}}dx\nonumber\\
&=&\lambda^{r-cp}\sum_{p=0}^{\infty}(-1)^{p}\binom{k+p-1}{p}\int_{\left(\frac{t}{\lambda}\right)^\beta}^{\infty}u^{1+\frac{r}{\beta}-\frac{cp}{\beta}-1}e^{-u}du\nonumber\\
&=& \lambda^{r-cp}\sum_{p=0}^{\infty}(-1)^{p}\binom{k+p-1}{p}\Gamma\left(\frac{r-cp}{\beta}+1,\left(\frac{t}{\lambda}\right)^{\beta}\right)
\end{eqnarray}
Substituting (\ref{c1}), (\ref{c2}), (\ref{c3}) in (\ref c), we get,
\begin{eqnarray}
E\left({X^{r}}/{X>t}\right)&=&\frac{1}{\bar{F(t)}}\left(k\sum_{m=0}^{\infty}\frac{(-1)^{m}}{\lambda^{m\beta}m!}B_{(1+t^{-c})^{-1}}\left(k+\frac{r}{c}+\frac{m\beta}{c},1-\frac{r}{c}-\frac{m\beta}{c}\right)\right)\nonumber
\\&+&\frac{1}{\bar{F(t)}}\left( \lambda^{r}\Gamma\left(\left(\frac{r}{\beta}+1\right),\left(\frac{t}{\lambda}\right)^{\beta}\right)\right)\nonumber
\\&-&\frac{1}{\bar{F(t)}}\left( \sum_{p=0}^{\infty}(-1)^{p}\binom{k+p-1}{p}\lambda^{r-cp}\Gamma\left(\frac{r-cp}{\beta}+1,\left(\frac{t}{\lambda}\right)^{\beta}\right)\right)
\end{eqnarray}

\vskip.3cm {\subsection*{3.4\hskip.3cm Mean Deviation}}\vskip.3cm

The amount of scatter in a population is measured to some extent by the totality of deviations from the mean and median. These are known as mean deviation about mean and as mean deviation about median and are defined as,
\begin{eqnarray*} \delta_1(x)&=&\int_{0}^{\infty} \mid{x-\mu}\mid f_{BW}(x)dx \ \ \text{and} \ \  \delta_2(x)=\int_{0}^{\infty} \mid{x-M}\mid f_{BW}(x)dx\end{eqnarray*} \\respectively where $ \mu=E(X)$ and \text{M$=$Median(X)} denote the median.
 The measures $\delta_1(x)$ and $\delta_2(x)$ can be calculated using the relationships,
  \begin{eqnarray}
 	\delta_1(x) &=&
 	2\mu F_{BW}(\mu)-2\mu+2\int_{\mu}^{\infty}xf_{BW}(x)dx
 	\\\delta_2(x)&=&-\mu+2\int_{M}^{\infty}xf_{BW}(x)dx
 	\end{eqnarray}respectively. When $r=1$ we get the mean $\mu=E(X)$. Note that $T(\mu)=\int_{\mu}^{\infty}xf_{BW}(x)dx$ and $T(M)=\int_{M}^{\infty}xf_{BW}(x)dx$, where
\begin{eqnarray}
T(\mu) &=& \int_{\mu}^{\infty}xf(x)dx\nonumber
\\&=&k\sum_{m=0}^{\infty}\frac{(-1)^{m}}{\lambda^{m\beta}m!}B_{(1+\mu^{-c})^{-1}}\left(k+\frac{1}{c}+\frac{m\beta}{c},1-\frac{1}{c}-\frac{m\beta}{c}\right)+\lambda\Gamma\left(\frac{1}{\beta}+1,\left(\frac{\mu}{\lambda}\right)^{\beta}\right)\nonumber
\\&-&\sum_{p=0}^{\infty}\lambda^{1-cp}(-1)^{p}\binom{k+p-1}{p}\Gamma\left(\frac{1-cp}{\beta}+1,\left(\frac{\mu}{\lambda}\right)^{\beta}\right)
\end{eqnarray}similarly,
\begin{eqnarray}
T(M) &=& \int_{M}^{\infty}xf(x)dx\nonumber
\\
\nonumber &=& k\sum_{m=0}^{\infty}\frac{(-1)^{m}}{\lambda^{m\beta}m!}B_{(1+M^{-c})^{-1}}\left(k+\frac{1}{c}+\frac{m\beta}{c},1-\frac{1}{c}-\frac{m\beta}{c}\right)\nonumber\\  &+&\lambda\Gamma\left(\frac{1}{\beta}+1,\left(\frac{M}{\lambda}\right)^{\beta}\right)-\sum_{p=0}^{\infty}
\lambda^{1-cp}(-1)^{p}\binom{k+p-1}{p}\Gamma\left(\frac{1-cp}{\beta}+1,\left(\frac{M}{\lambda}\right)^{\beta}\right)\nonumber\\
\end{eqnarray}
Consequently, the mean deviation about the mean and the mean deviation about the median reduces to \begin{eqnarray*}\delta_1(x)&=&2\mu F_{BW}(\mu)-2\mu+2T(\mu) \ \ \ and \ \ \delta_2(x)=-\mu+2T(M)\end{eqnarray*}respectively.

\vskip.3cm {\subsection*{3.5\hskip.3cm Bonferroni and Lorenz curves}}\vskip.3cm

 Bonferroni and Lorenz curves have applications not only in economics for the study of income and poverty, but also in other fields such as reliability, demography, insurance and medicine. Bonferroni and Lorenz curves for the BW distribution are given by,
 \begin{eqnarray*}
 	B(p)&=& \frac{1}{p\mu}\int_{0}^{q}xf_{BW}(x)dx =\frac{1}{p\mu}[\mu-T(q)],\ \ \ \text{and} \\L(p) &=&\frac{1}{\mu}\int_{0}^{q}xf_{BW}(x)dx=\frac{1}{\mu}[\mu-T(q)],\end{eqnarray*}\\respectively, where
 \begin{eqnarray}
 T(q)&=&\int_{q}^{\infty}xf_{BW}(x)dx\nonumber\\
&=&k\sum_{m=0}^{\infty}\frac{(-1)^{m}}{\lambda^{m\beta}m!}B_{(1+q^{-c})^{-1}}\left(k+\frac{1}{c}+\frac{m\beta}{c},1-\frac{1}{c}-\frac{m\beta}{c}\right)\nonumber\\ &+&\lambda\Gamma\left(\frac{1}{\beta}+1,\left(\frac{q}{\lambda}\right)^{\beta}\right)-\sum_{p=0}^{\infty}\lambda^{1-cp}(-1)^{p}\binom{k+p-1}{p}
\Gamma\left(\frac{1-cp}{\beta}+1,\left(\frac{q}{\lambda}\right)^{\beta}\right)
\end{eqnarray}and $q=F^{-1}(p),0\leq p\leq 1$.

\vskip.3cm {\subsection*{3.6\hskip.3cm Order Statistics}}\vskip.3cm

The density function $f_{i:m}(x)$ of the $i^{th}$ order statistic for i=1,2,3,...,m from independently and identically distributed random variables
$X_{1},X_{2},...,X_{m}$ following BW distribution is given by,
\begin{eqnarray}\label{os}
f_{i:m}(x) &=& \frac{m!f_{BW}(x)}{(i-1)!(m-i)!}[F_{BW}(x)]^{i-1}[1-F_{BW}(x)]^{m-i}\nonumber\\
&=&\frac{m!f_{BW}(x)}{(i-1)!(m-i)!}\sum _{j=0}^{m-i}(-1)^{j}\binom{m-i}{j}[F_{BW}(x)]^{j+i-1}
\end{eqnarray}
Using the binomial expansion $[1-F(x)]^{m-i}=\sum_{j=0}^{m-i}\binom{m-i}{j}(-1)^{j}[F(x)]^{j}$
and the pdf and cdf of BW distribution in (\ref{os}) we have,
\begin{eqnarray*}
f_{i:m}(x) &=&\frac{m!f_{BW}(x)}{(i-1)!(m-i)!}\sum _{j=0}^{m-i}(-1)^{j}\binom{m-i}{j}\left(1-\left[1-\left(1+x^{-c}\right)^{-k}\right]\left[e^{-\left(\frac{x}{\lambda}\right)^{\beta}}\right]\right)^{j+i-1}\\
&=&\sum_{j=0}^{m-i}(-1)^{j}\frac{m!}{(i-1)!(m-i-j)!(j)!}\left(1-\left[1-\left(1+x^{-c}\right)^{-k}\right]\left[e^{-\left(\frac{x}{\lambda}\right)^{\beta}}\right]\right)^{j+i-1}\\
&\times&e^{-\left(\frac{x}{\lambda}\right)^{\beta}}
\left(ck\left(1+x^{-c}\right)^{-k-1}x^{-c-1}+\frac{\beta}{\lambda^{\beta}}x^{\left(\beta-1\right)}\left(1-\left(1+x^{-c}\right)^{-k}\right)\right)
\end{eqnarray*}
The pdf of the $1^{st}$ and $n^{th}$ order statistic is given by,
\begin{eqnarray}
f_{1:m}(x)&=&\sum_{j=0}^{m-1}(-1)^{j}\frac{m!}{(m-1-j)!(j)!}\left(1-\left[1-\left(1+x^{-c}\right)^{-k}\right]\left[e^{-\left(\frac{x}{\lambda}\right)^{\beta}}\right]\right)^{j}\nonumber\\
&\times&e^{-\left(\frac{x}{\lambda}\right)^{\beta}}
\left(ck\left(1+x^{-c}\right)^{-k-1}x^{-c-1}+\frac{\beta}{\lambda^{\beta}}x^{\left(\beta-1\right)}\left(1-\left(1+x^{-c}\right)^{-k}\right)\right)\\
f_{n:m}(x)&=&\sum_{j=0}^{m-n}(-1)^{j}\frac{m!}{(n-1)!(m-n-j)!(j)!}\left(1-\left[1-\left(1+x^{-c}\right)^{-k}\right]\left[e^{-\left(\frac{x}{\lambda}\right)^{\beta}}\right]\right)^{j+n-1}\nonumber\\
&\times&e^{-\left(\frac{x}{\lambda}\right)^{\beta}}
\left(ck\left(1+x^{-c}\right)^{-k-1}x^{-c-1}+\frac{\beta}{\lambda^{\beta}}x^{\left(\beta-1\right)}\left(1-\left(1+x^{-c}\right)^{-k}\right)\right)
\end{eqnarray}

\vskip.5cm{\section*{4.\hskip.3cm Renyi's Entropy}}
\vskip.3cm
In this section, Renyi's entropy of the BW distribution is derived. An entropy is a measure of uncertainty or disorder of a random variable. Renyi's entropy is an extension of Shannon's entropy. In the case of BW distribution Renyi's entropy is defined to be
 \begin{eqnarray*}
	I_R(v)&=&\frac{1}{1-v}\ln\left(\int_{0}^{\infty}[f_{BW}(x;c,k,\lambda,\beta)]^v dx \right),v\neq1, v>0.\end{eqnarray*}
Renyi's entropy tends to Shannon's entropy as $v\rightarrow 1$. Note that $[f(x;c,k,\lambda,\beta)]^v=f_{BW}^v(x)$ can be written as,\begin{eqnarray*} f_{BW}^v(x)&=&e^{-v\left(x/\lambda\right)^{\beta}} \left[ck\left(1+x^{-c}\right)^{-k-1}x^{-c-1}+\frac{\beta}{\lambda^\beta}x^{\beta-1}\left(1-\left(1+x^{-c}\right)^{-k}\right)\right]^{v}
	 \end{eqnarray*}
	 Using the expansions,\begin{eqnarray*} e^{-x}&=&\sum_{j=0}^{\infty}\frac{(-1)^j x^j}{j!}, \ (x+y)^n=\sum_{p=0}^{n} \  ^nC_p x^{n-p}y^p, \ (1-x)^p=\sum_{m=0}^{p}  \  ^pC_m (-1)^m x^m. \end{eqnarray*}
\begin{eqnarray*} f_{BW}^v(x)&=&\sum_{j=0}^{\infty} \frac{(-1)^jv^j (\frac{x}{\lambda})^{\beta j}}{j!}\left[\frac{\beta x^{\beta-1}}{\lambda^\beta}(1-(1+x^{-c})^{-k})+kcx^{-c-1}(1+x^{-c})^{-k-1}\right]^v
	 	\\&=&\sum_{j=0}^{\infty} \frac{(-1)^j v^jx^{\beta j}}{\lambda^{\beta j}j!} \sum_{p=0}^{v} \ {{v}\choose{p}} \left(kcx^{-c-1}(1+x^{-c})^{-k-1}\right)^{v-p} \left(\frac{\beta x^{\beta-1}}{\lambda^\beta}(1-(1+x^{-c})^{-k})\right)^p
	 	 \\&=& \sum_{j=0}^{\infty}\sum_{p=0}^{v}\sum_{w=0}^{p}\frac{(-1)^{j+w}  v^j  \beta^ p}{\lambda^{(\beta j+p\beta)} j!}(kc)^{v-p} {{v}\choose{p}}  {{p}\choose{w}} x^{\beta j+cp+p-cv-v+p\beta-p} \\&\times&(1+x^{-c})^{kp+p-kv-v-kw}  .  \end{eqnarray*}
	 	\\Now, \begin{eqnarray*}\int_{0}^{\infty}f_{BW}^v(x)dx&=&\sum_{j=0}^{\infty}\sum_{p=0}^{v}\sum_{w=0}^{p}\frac{(-1)^{j+w} v^j \beta^ p}{\lambda^{(\beta j+p\beta)} j!}(kc)^{v-p} {{v}\choose{p}}  {{p}\choose{w}}\\&\times& \int_{0}^{\infty}x^{\beta j+cp+p-cv-v+p\beta-p} (1+x^{-c})^{kp+p-kv-v-kw} dx. \end{eqnarray*}
	 	Put $u=(1+x^{-c})^{-1}$ , \begin{eqnarray*}
	 		\int_{0}^{\infty}f_{BW}^v(x)dx&=&\sum_{j=0}^{\infty}\sum_{p=0}^{v}\sum_{w=0}^{p}\frac{(-1)^{j+w} v^j \  \beta^ p}{\lambda^{(\beta j+p\beta)} \ j!} \ (kc)^{v-p} \ {{v}\choose{p}} \ {{p}\choose{w}}\ \frac{1}{c}\\&\times& \int_{0}^{\infty}u^{a-1} (1-u)^{b-1} du \\&=&\sum_{j=0}^{\infty}\sum_{p=0}^{v}\sum_{w=0}^{p}\frac{(-1)^{j+w}\ v^j \ \beta^ p}{\lambda^{(\beta j+p\beta)}\ j!}\ (kc)^{v-p}\ {{v}\choose{p}}\  {{p}\choose{w}} \ \frac{1}{c}\times B(a,b),
	 		 \end{eqnarray*}where
	 		 \begin{eqnarray*}
	 		a&=&\frac{1}{c}(p\beta+\beta j+cp-cv-v)+\frac{1}{c}-kp-p+kv+v+kw, \\b&=&\frac{1}{c}(cv+v-p\beta-\beta j-cp)-\frac{1}{c} .\end{eqnarray*}
	 	Then,
	 	\begin{eqnarray}
	 		I_R(v)&=&\frac{1}{1-v}\ln \left(\sum_{j=0}^{\infty}\sum_{p=0}^{v}\sum_{w=0}^{p}\frac{(-1)^{j+w}\ v^j \ \beta^ p}{\lambda^{(\beta j+p\beta)}\ j!}\ (kc)^{v-p}\ {{v}\choose{p}}\  {{p}\choose{w}} \ \frac{1}{c} B(a,b)\right),
	 		\end{eqnarray} \\for $v\neq 1$ abd $v>0$.

\vskip.5cm{\section*{5.\hskip.3cm Method of Maximum Likelihood Estimation}}
\vskip.3cm

The most useful parametric estimation method is the maximum likelihood method. Let $x_{1},x_{2},....x_{n}$ be a random sample of size n from BW distribution. Then the  log likelihood function is given by,
\begin{eqnarray}
l(c,k,\lambda,\beta)&=&\sum_{i=1}^{n}\ln [f(x_{i},c,k,\lambda,\beta)]
\end{eqnarray}
The log likelihood for a single observation is given by,
\begin{eqnarray}\label{mle}
l(c,k,\lambda,\beta)&=&-\left(\frac{x}{\lambda}\right)^{\beta}+\ln(ck(1+x^{-c})^{-k-1}x^{-c-1}+\frac{\beta}{\lambda^\beta}x^{\beta-1}[1-(1+x^{-c})^{-k}])
\end{eqnarray}
The maximum likelihood estimates can be obtained by solving the following the equations simultaneously,\\$\frac{\partial l(c,k,\lambda,\beta)}{\partial c}=0,\frac{\partial l(c,k,\lambda,\beta)}{\partial k}=0,\frac{\partial l(c,k,\lambda,\beta)}{\partial \lambda}=0,\frac{\partial l(c,k,\lambda,\beta)}{\partial \beta}=0$ where
\begin{align*}
\resizebox{0.9\hsize}{22pt}	{$\frac{\partial l(c,k,\lambda,\beta)}{\partial c} =\frac{k\lambda^{\beta}\left(-cx^{-c-1}\ln x\left(1+x^{-c}\right)^{-k-1}-cx^{-2c-1}(-k-1)(1+x^{-c})^{-k-2}\ln x+x^{-c-1}(1+x^{-c})^{-k-1}\right)-k\beta x^{-c+\beta-1}(1+x^{-c})^{-k-1}\ln x}{\lambda^{\beta}ckx^{-c-1}(1+x^{-c})^{-k-1}+\frac{\beta}
{\lambda^\beta}x^{\beta-1}\left[1-\left(1+x^{-c}\right)^{-k}\right]}$}
\end{align*}

\begin{align*}
\resizebox{0.9\hsize}{22pt}{$\frac{\partial l(c,k,\lambda,\beta)}{\partial k }=\frac{ckx^{-c-1}(1+x^{-c})^{-k-1}-cx^{-c-1}(1+x^{-c})^{-k-1}\ln(1+x^{-c})+\frac{\beta}{\lambda^\beta}x^{\beta-1}\ln(1+x^{-c})(1+x^{-c})^{-k}}{ckx^{-c-1}(1+x^{-c})^{-k-1}+\frac{\beta}{\lambda^\beta}x^{\beta-1}\left[1-\left(1+x^{-c}\right)^{-k}\right]}$}
\end{align*}
\begin{align*}
\resizebox{0.9\hsize}{22pt}	{$\frac{\partial l(c,k,\lambda,\beta)}{\partial \lambda }=\frac{\beta x}{\lambda^{2}}\left(\frac{x}{\lambda}\right)^{\beta-1}-\frac{\beta ^{2}x^{\beta-1}\lambda^{-\beta-1}\left[1-\left(1+x^{-c}\right)^{-k}\right]} {ckx^{-c-1}(1+x^{-c})^{-k-1}+\frac{\beta}{\lambda^\beta}x^{\beta-1}\left[1-\left(1+x^{-c}\right)^{-k}\right]}$}
\end{align*}
\begin{align*}
\resizebox{0.9\hsize}{22pt}	{$\frac{\partial l(c,k,\lambda,\beta)}{\partial \beta }= -\left(\frac{x}{\lambda}\right)^{\beta}\ln\left(\frac{x}{\lambda}\right)+\frac{\left[1-\left(1+x^{-c}\right)\right]^{-k}\left[\lambda^{2}\left(x^{\beta-1}+\beta x^{\beta-1}\text{log(x)}(x)\right)-\beta x^{\beta-1} \lambda^\beta \ln\lambda\right]}{(\lambda^{\beta})^{2}\left[ckx^{-c-1}(1+x^{-c})^{-k-1}+\frac{\beta}{\lambda^\beta}x^{\beta-1}\left[1-\left(1+x^{-c}\right)^{-k}\right]\right]}$}
\end{align*}
The total log likelihood function based on random sample of n observations $x_{1},x_{2},....x_{n}$ drawn from BW distribution is given by $l^{*}=\sum_{i=1}^{n}l(c,k,\lambda,\beta)$ where $l(c,k,\lambda,\beta)$ is given by equation (\ref{mle}). Owing to the complexity of these equations, the MLEs does not have an analytical expression. However, one can use standard statistical software to solve those equations (e.g., Mathematica, R, etc.). We make use of R software to carry out this study.

\vskip.5cm{\section*{6.\hskip.3cm Simulation}}
\vskip.3cm
{\scriptsize{
\begin{center}
\begin{table}
\caption{Simulation Results: Mean Estimates, Bias, MSE.}
	\label{tt}
	\begin{tabular}{|c|c| c| c| c|c|c| c |c|}
		\hline
		&&\multicolumn{3}{|c|}{I}&	\multicolumn{3}{c|}{II} \\
		\hline
		&&&&&&&\\
		Sample Size & Parameter & Mean& Bias& MSE & Mean & Bias & MSE \\
\hline	
		\multirow{4}{1pt}{n=25}
		& c &0.3843&0.08613& 0.02274 & 5.9395 &0.43952& 1.79489\\
		&k & 7.8764 &-0.10006 &0.22994&5.0623  &0.06233 &0.68445 \\
		&$\lambda$ &1.1783 &-0.03159&0.01312&0.89611 & -0.00388& 0.00353\\
		&$\beta$ & 2.32138&0.297792&0.22436 & 3.6348&  0.33482& 0.47132 \\
		\hline
		
		\multirow{4}{1pt}{n=200}
		& c & 0.3557&0.05578& 0.00538 & 5.7685 &0.2685& 0.7198\\
		&k & 7.9368 &-0.06318 &0.06149&5.06035  & 0.06035 &0.12072\\
		&$\lambda$ & 1.1776 & -0.02231 &0.00211&0.8961 & -0.0038& 0.0003\\
		&$\beta$ & 2.2158&0.21584 &0.05967 & 3.5156&  0.2156& 0.0854 \\
		\hline
		\multirow{4}{1pt}{n=400}
		&c&0.3538&0.05388&0.004074&5.7476 &0.2476& 0.2950\\
		&k&  7.9723&-0.02077&0.003812&5.0564  &0.05644 &  0.03913\\
		&$\lambda$ &1.1769&-0.02184&0.00129 &0.8964  & -0.0035&0.00019\\
		&$\beta$&2.2105 &0.210585& 0.050837& 3.5054&  0.20548&0.05793\\
		\hline
		\multirow{4}{1pt}{n=600}
		&c& 0.3534& 0.05343&0.003653& 5.6836&0.1836& 0.1939\\
		&k& 7.9823&-0.01765&0.001241&5.0501  &0.05018 &  0.03712\\
		&$\lambda$ &1.1767&-0.02158&0.00098 &0.8963 & -0.0034&0.00014\\
		&$\beta$&2.2057&0.20573& 0.046781& 3.5039&  0.2039&0.0533\\
		\hline
			\end{tabular}
\end{table}
\end{center}}}
The accuracy and performance of the BW distribution is investigated by conducting two
simulations for different parameter values and sample sizes. The simulations were repeated $N = 1000$
times each with sample sizes $n=25, 50, 200, 400, 600$ and the true
parameters values $\text{I}:c=0.3,k=8,\lambda=1.2,\beta=2$ and $ \text{II}:c=5.5,k=5,\lambda=0.9,\beta=3.3$.
Three quantities were computed in this simulation study: the mean, bias and mean-square error (MSE).
The mean estimate, bias and mean-square error of the MLE $\hat{\epsilon}$ of the parameter $c, k, \lambda, \beta$ are respectively  given by


$$
\text{Mean} =\frac{\sum_{i=1}^{N} (\hat{\epsilon_i})}{N}, ~\text{Bias} =\frac{1}{N}\sum_{i=1}^{N}(\hat{\epsilon}-\epsilon) \quad \text{and} \quad \text{MSE}=\frac{1}{N}\sum_{i=1}^{N}(\hat{\epsilon}-\epsilon)^{2}.
$$

The mean of MLEs of the BW distribution parameters along with their respective
mean square errors and bias for different sample sizes are listed in Tables \ref{tt}.

%
\vskip.5cm{\section*{7.\hskip.3cm Real Applications}}
\vskip.3cm

This section illustrates the usefulness of the Burr III Weibull distribution using a real data set.
The data set includes 101 observations
which represent the failure times of Kevlar 49/epoxy
strands which were subjected to constant sustained pressure
at the 90 \% percent stress level until all had failed. The data are: 0.01, 0.01, 0.02, 0.02, 0.02,0.03, 0.03, 0.04, 0.05, 0.06, 0.07, 0.07, 0.08, 0.09, 0.09, 0.10, 0.10, 0.11, 0.11, 0.12, 0.13, 0.18,0.19, 0.20, 0.23, 0.24, 0.24, 0.29, 0.34, 0.35, 0.36, 0.38, 0.40, 0.42, 0.43, 0.52, 0.54, 0.56, 0.60,0.60, 0.63, 0.65, 0.67, 0.68, 0.72, 0.72, 0.72, 0.73, 0.79, 0.79, 0.80, 0.80, 0.83, 0.85, 0.90, 0.92,0.95, 0.99, 1.00, 1.01, 1.02, 1.03, 1.05, 1.10, 1.10, 1.11, 1.15, 1.18,1.20, 1.29, 1.31, 1.33, 1.34,1.40, 1.43, 1.45, 1.50, 1.51, 1.52, 1.53, 1.54, 1.54, 1.55, 1.58, 1.60, 1.63, 1.64, 1.80, 1.80, 1.81,2.02, 2.05, 2.14, 2.17, 2.33, 3.03, 3.03, 3.34, 4.20, 4.69, 7.89.

\vskip.3cm We have fitted the proposed  Burr III Weibull distribution (BW) to the data set and compared the proposed distribution with Weibull distribution and Burr III distribution. The values of the estimated parameters, Akaike Information Criterion (AIC), Bayesian Information Criterion (BIC) and Consistent Akaike Information Criterion (AICC) values for the corresponding data are provided in Table \ref{Table 2}. The BW distribution is a better model as compared to the Burr III model.  We conclude that the Burr III-Weibull distribution can be comparable to the Burr III and Weibull models.

{\scriptsize{
\begin{table}
	\caption{Comparison of Maximum Likelihood Estimates.}
	\label{Table 2}
	\begin{tabular}{|c|c|c|c|c|c|c|}
		\hline
		\multirow{2}{*}{Model}
		&\multicolumn{2}{c|}{MLE}&\multirow{2}{*}{Log-Liklyhood}&\multirow{2}{*}{AIC}&\multirow{2}{*}{BIC}&\multirow{2}{*}{AICC}\\
		&Parameters&Estimates&&&&\\
		\hline
		BurrIII&c&2.3858890&-98.66771&205.3354&215.7959&205.7521\\
		-Weibull&k&2.4533820&&&&\\
		&$\lambda$&1.7572900&&&&\\
		&$\beta$&0.6791234&&&&\\
		\hline
		BurrIII&c&1.8321566&-106.6097&217.2194&222.4497&217.3419\\
		&k&0.5343506&&&&\\
		\hline
		Weibull&$\lambda$&0.9899448&-102.9768&209.9536&215.1839& 210.0761\\
		&$\beta$&0.9258876&&&&\\
		\hline
	\end{tabular}
\end{table}}}

\vskip.5cm{\section*{8.\hskip.3cm Conclusion}}
\vskip.3cm

A new distribution called Burr III-Weibull(BW)distribution  is proposed and its properties are
studied. The BW distribution possesses increasing, decreasing and upside-down bathtub shaped failure rate function. We derived the moments, conditional moments, mean deviation, quantiles, Bonferroni and Lorentz curve etc. of the proposed distribution. Order statistics and Renyi's entropy of the proposed distribution are  also obtained. Estimation of the parameters of the distribution is performed via maximum likelihood method.  A simulation study is performed to validate the maximum likelihood estimator (MLE). Finally, the BW distribution
is fitted to real data sets in order to illustrate its applicability and usefulness.

\end{document}